\documentclass[12pt]{amsart}
\usepackage{amsmath}\usepackage{amsrefs}
\usepackage{amsfonts}\usepackage{verbatim}
\usepackage{amssymb}\usepackage{marvosym}\usepackage{skak}\usepackage{dictsym}
\usepackage{amstext}
\usepackage{amsbsy}
\usepackage{amsopn}\usepackage{MnSymbol}
\usepackage{amsthm}

\theoremstyle{definition}

\theoremstyle{remark}

\DeclareMathOperator\const{const}
\def\proclaim#1{\vskip0.5em\noindent{\bf #1}\it}
\def\endproclaim{\vskip0.5em\par\noindent\rm}
\def\proclaim#1{\vskip0.5em\noindent{\bf #1}\it}
\def\endproclaim{\vskip0.5em\par\noindent\rm}
\def\demo#1{\vskip0.5em\noindent{\bf #1\ }}

\def\cite#1{[#1]}
\def\text#1{\mbox{#1}}
\def\flushpar{\par\noindent}

\newcommand{\mapright}[1]{%
    \smash{\mathop{%
        \hbox to 1cm{\rightarrowfill}
        }
    \limits^{#1}
    }
}
\newcommand{\mapleft}[1]{%
    \smash{\mathop{%
        \hbox to 1cm{\rightarrowfill}
        }
    \limits_{#1}
    }
}

\def\e{\epsilon}
\def\a{\alpha}

\def\G{\Gamma}
\def\g{\gamma}
\def\d{\delta}

\def\s{\sigma}
\def\Si{\Sigma}
\def\th{\theta}
\def\l{\lambda}
\def\x{\times}

\def\o{\overline}
\def\f{\flushpar}
\def\un{\underline}
\def\vf{\varphi}

\def\om{\omega}
\def\Om{\Omega}
\def\B{\mathcal B}

\def\({\biggl(}
\def\){\biggr)}\def\bdy{\partial}

\def\<{\langle}
\def\>{\rangle}
\def\bul{\smallskip\f$\bullet\ \ \ $}

\def\lfl{\lfloor}\def\rfl{\rfloor}\def\sms{\smallskip\f}
\def\sms{\smallskip\f}\def\sbul{\f$\bullet\
\ \ $}\def\sms{\smallskip\f}
\def\lra{\longrightarrow}\def\Lra{\Longrightarrow}\def\lfl{\lfloor}\def\rfl{\rfloor}\def\st{\text{such that}}
\begin{document} \title{ Exponential chi squared distributions in infinite ergodic theory}
\author{ Jon. Aaronson $\&$ Omri Sarig}\address[Jon. Aaronson]{\ \ School of Math. Sciences, Tel Aviv University,
69978 Tel Aviv, Israel.}
\email{aaro@tau.ac.il}
\address[Omri Sarig]{Faculty of Mathematics and Computer Sciences,
The Weizmann Institute for Science, POB 26, Rehovot 76100, Israel }
\email{omsarig@gmail.com}

\begin{abstract} We prove  distributional limit theorems  for  random walk adic transformations obtaining   ergodic distributional limits of  exponential chi squared form.
\end{abstract}

\thanks{Aaronson's research was
supported by Israel Science Foundation grant No. 1114/08.\ Sarig's  research  was supported by the European Research Council, grant 239885}
\maketitle\markboth{ Jon. Aaronson $\&$ Omri Sarig \copyright  2012}{Exponential $\chi^2$ distributions}

\section*{\S0 Introduction}
\

As in [A1], for $(X,\B,m)$  a $\s$-finite measure space, $F_n:X\to [0,\infty]$ measurable, and $Y\in [0,\infty]$ a random variable,
we say that {\it $(F_n)$ converges strongly in distribution to $Y$}, (written $F_n\overset{\mathfrak d}{\underset{n\to\infty}\lra } Y ,$)
if it converges in law with respect to all $m$-absolutely continuous probabilities; equivalently
$$g(F_n)\underset{n\to\infty}\lra \, \Bbb E(g(Y))\ \text{\rm weak $*$ in}\ L^\infty(m)$$ for each bounded, continuous function $f:[0,\infty]\to\Bbb R$.
\

For discussion of strong distributional convergence, see [A1], [A2], [E] and [TZ].
\

Here, we study {\tt distributional stability}. As in [A2], we'll call a conservative, ergodic measure preserving transformations
 $(X,\B,m,T)$   {\it  distributionally stable} if there are constants $a(n)>0$ and a random variable $Y$ on $(0,\infty)$ so that
$$\tfrac1{a(n)}S_n(f)\overset{\mathfrak d}\lra \ Ym(f)\ \forall\ f\in L^1_+$$  where
$S_n(f):= \sum_{k=0}^{n-1} f\circ T^k$ and $m(f):=\int_Xfdm$.
\par By the ratio ergodic theorem, if the above convergence holds for some $f\in L^1_+$, then it holds $\forall\ f\in L^1_+$.
\par If the {\it ergodic distributional limit} $Y$  is normalized by $\Bbb E(Y)=1$, the constants  $a(n)$ are determined uniquely up to asymptotic equality and are aka the {\it return sequence} of $T$. Both the (normalized) ergodic distributional limit and the return sequence  are invariant under similarity  (see [A1]).
\

By the Darling-Kac theorem ([DK]), pointwise dual ergodic transformations (e.g. Markov shifts) with regularly varying return sequences are distributionally stable with Mittag-Leffler  ergodic distributional limits (see also [A1], [A2]).

\

Our present study concerns {\tt\small random walk adic } {\tt\small transformations}.

\

A {\it random walk adic transformation} is a conservative, ergodic measure preserving transformation  associated to a Markov driven, aperiodic,  random walk on a group of form $\Bbb G=\Bbb Z^k\x\Bbb R^{D-k}$.
\

It is  the (unique) $\Bbb G$-extension of the adic transformation on the underlying Markov shift which parametrizes the tail relation of the random walk. This definition is explained in \S1.
\

The  {\it degree} of a random walk adic transformation is the
dimension of the associated group $\Bbb G$: $\dim(\text{\tt span}\,\Bbb G)$. It appears as the number of degrees of freedom in the $\chi^2$ distribution appearing in the limit.
\

 The  ``original HIK transformation" appears in [HIK].
\

We establish the following
\proclaim{Theorem}
\

Suppose that $(X,\B,m,T)$ is a random walk adic transformation with degree $D\in\Bbb N$,
then
$$\tfrac1{a_n(T)}S_n(f)\overset{\mathfrak d}\lra \ (2^{\frac{D}{2}}e^{-\frac{1}{2}\chi^2_{D}})\mu(f)\ \forall\ f\in L^1_+$$ where
$S_n(f):= \sum_{k=0}^{n-1} f\circ T^k$, $a_n(T)\propto\frac{n}{(\log n)^{D/2}}$; and $\chi^2_{D}=\|\xi\|_2^2$ for $\xi$ a standard Gaussian random vector on $\Bbb R^D$.\endproclaim

This ergodic distributional limit first appears in [LS] (see below).
\

Most of the paper is devoted to proving the theorem. In \S1, we define adic transformations and  random walk adic transformations as group extensions of adic transformations. In  \S2, we establish compact representation for (all) adic transformations and a uniform convergence theorem for stationary adic transformations which latter is needed in the proof of the theorem. We review the distributional limit theory of Markov shifts in \S3 and prove the theorem in \S4, giving  applications to exchangeability in \S5 and horocycle flows in \S6.

\

Related, earlier work can be found in [AW], [ANSS] and  [LS] (see \S6).
\section*{\S1 Bratteli diagrams, adic  and  random walk adic transformations}
\

\subsection*{Bratteli diagrams}
\

Fix $a_n\ge 2\ \ (n\ge 1)$ and set $\mathcal S_k:=\{0,1,\dots,a_k-1\},\  \Om:=\prod_{k=1}^\infty\mathcal S_k$.
\

A {\it Bratteli diagram} is a subset $\Sigma\subset\Om$ of form
$$\Sigma:=\{\om\in\Om:\ \ A_k(\om_k,\om_{k+1})=1\ \ \forall\ k\ge 1\}$$
where for $k\ge 1$,
$A_k:\mathcal S_k\x\mathcal S_{k+1}\to\{0,1\}$ is the {\it $k^{\text{\tiny th}}$ transition matrix}.
\

Recall that $\Om$ is compact when equipped with the standard metric $d(x,y):=\exp[-\min\{n:x_n\neq y_n\}]$ and $\Si$ is a closed subset.
\

The Bratteli diagram $\Si$ is called {\it stationary} if  $\mathcal S_k=\mathcal S,\ A_k=A\ \forall\ k\ge 1$. In this case,
 $\Si$ is a {\it topological Markov shift} (TMS) with {\it transtition matrix} $A$ as in [LM].
 \

 The only result in this paper concerning non-stationary Bratteli diagrams is the compact representation lemma in \S2.
\

\subsection*{Tail relation on a Bratteli diagram}

\

The {\it tail relation} on $\Si$  is the equivalence relation
$$\mathfrak T=\mathfrak T(\Si):=\{(x,y)\in\Si\x\Si:\ \exists n\ \ \st\ x_n^\infty=y_n^\infty\}$$ where
$x_n^\infty:=(x_n,x_{n+1},\dots)$.
\

The equivalence classes of the tail relation are linearly ordered by the {\it reverse lexicographic order}, namely
the  partial order $\prec$ on $\Sigma$ defined  by
$$
x\prec y\Leftrightarrow \exists n\textrm{ s.t. }x_{n+1}^{\infty}=y_{n+1}^{\infty}\textrm{ and }x_n< y_n.
$$
If  $x$ is maximal, then $x_n=\max\{y\in \mathcal S_n:A_n(y,x_{n+1})=1\}\ \forall\ n\ge 1$; therefore the collection of non-maximal points is open and the collection $\Si_{\text{\tt\tiny max}}$ of maximal points is closed.
 A similar argument shows that the set $\Si_{\text{\tt\tiny min}}$ of minimal points is closed.
\

In case $\Si$ is a  TMS (stationary Bratteli diagram), more is true.
\

If  $x$ is maximal, then $x_n=\vf_+(x_{n+1})$ where $\vf_+(x)=\max\{y\in \mathcal S:A(y,x)=1\}$; and we claim that \begin{align*}\tag{{\Large\dsrailways}}\ \ \ \text{\large\emph{ $x$ is periodic, with period $\le\#\mathcal S$.}}\end{align*}
\

To see ({\Large\dsrailways}), note first that $\exists\ s\in\mathcal S$ so that $\#\{n\ge 1:\ x_n=s\}=\infty$. The sequence $n\mapsto\vf_+^n(s)$ is eventually periodic with a final period $$(t,\vf_+(t)\dots,\vf_+^{\kappa-1}(t))=(\vf_+^{J}(s),\dots,\vf_+^{J+\kappa-1}(s))$$ with $J\ge 1\ \&\ \kappa\le \#\mathcal S$.
\

Let $p:=(\vf_+^{\kappa-1}(t),\vf_+^{\kappa-2}(t),\dots,t,\vf_+^{\kappa-1}(t),\vf_+^{\kappa-2}(t),\dots,t,\dots)$. We prove our  claim by showing that  $x=\s^\ell(p)$ for some $\ell\ge 1$. Let $n_k\uparrow\infty$ be so that $x_{n_k}=s\ \ \forall\ k\ge 1$. It follows that
$x_{n_k-\nu}=\vf_+^{\nu}(s)$ whence $\exists\ 1\le\ell_k\le\kappa$ so that
$x_1^{n_k-J}=\s^{\ell_k}(p)_1^{n_k-J}.$
There is a subsequence $m_j=n_{k_j}\to\infty$ so that $\ell_{k_j}=\ell_{k_1}=:\ell\ \forall\ j\ge 1$ whence
$x=\s^\ell(p)$ and ({\Large\dsrailways}) is established.
\

Thus the set of  maximal points (and the set of minimal points) is finite.
\subsection*{Adic transformations}
\

The {\it adic transformation} (generated by $\prec$) on the Bratteli diagram $\Si$  is $$\tau:\Si\setminus\Si_{\max}\to\Si\setminus\Si_{\min}\ \text{\tt defined by}\ \
\tau(x):=\min\{y:y\succ x\}.
$$
It is called {\it stationary} if the  underlying Bratteli diagram is stationary.
\

As shown in [V], every ergodic, probability preserving transformation is isomorphic to some adic transformation.
\

Stationary adic transformations are \sbul isomorphic to  odometers or primitive substitutions, have zero entropy but can be weakly mixing (see [L1]); and \sbul are always {\it uniquely ergodic}; moreover  the unique $\tau$--invariant probability measure $\nu_0$ is globally supported,   non-atomic, Markov and equivalent to the Parry  measure $\mu$ (of maximal entropy) for the associated TMS (see [BM]).

\subsection*{The $(\Si,f)$-random walk}
\

Let $\Sigma$ be a topologically mixing TMS on the (ordered) finite state space $\mathcal S$, let $\s:\Si\to\Si$ is the shift and let $\tau:\Si'\to\Si'$ be the corresponding stationary adic transformation
where
$$\Si':=\Si\setminus\bigcup_{n\ge 0}\s^{-n}(\Si_{\text{\tt\tiny max}}\cup\Si_{\text{\tt\tiny min}})=\bigcap_{n\in\Bbb Z}\tau^{n}(\Si\setminus(\Si_{\text{\tt\tiny max}}\cup\Si_{\text{\tt\tiny min}})).$$

The tail relation of $\Si$ is the tail relation of $\s$:
$$\mathfrak T(\Si)=\mathfrak T(\s):=\bigcup_{n\ge 0}\{(x,y)\in\Si\x\Si:\ \s^n(x)=\s^n(y)\};$$ and this is parametrized by the adic transformation:
 $$\mathfrak T(\s)\cap(\Si'\x\Si')=\{(x,\tau^n(x)):\ x\in\Si',\ n\in\Bbb Z\}.$$

\

A function $f:\Sigma\to\Bbb R^d$ is {\it H\"older continuous} if
 $\exists\ \th\in (0,1),\ M>1$ so that
 \begin{align*}\tag{\dsagricultural}\|f(x)-f(y)\|\le M\th^n\ \ \forall\ x,y\in\Si,\ x_n^\infty=y_n^\infty.\end{align*}
Specifically for we call $f:\Sigma\to\Bbb R^d$ {\it $\th$-H\"older continuous} ($\th\in (0,1)$) if (\dsagricultural) is satisfied for some $M>1$.
\

 For  $f:\Sigma\to\Bbb R^d$  H\"older continuous let
 $$\Bbb H:=\o{\<\{f_n(x):\ n\ge 1,\ x\in \Si,\ \s^nx=x\}\>}$$
 where
 $$f_n(x)=f^{(\s)}_n(x):=\sum_{k=0}^{n-1}f(\s^kx)$$
 and let
 $$\Bbb G:=\o{\<\{f_n(x)-f_n(y):\ n\ge 1,\ x,\ y\in \Si,\ \s^nx=x,\ \s^ny=y\}\>},$$
 then $\Bbb G,\ \Bbb H$ are both closed subgroups of $\Bbb R^d$ and $\Bbb G\le \Bbb H$.
 \

 It follows from  Liv\v{s}ic's cohomology theorem [L2], (see e.g. [ANS], [SA], [PS]) that $$f=g-g\circ\s+h+a\ \text{where}$$
  \bul $g:\Si\to\Bbb R^d$ is H\"older continuous;
  \bul $a\in\Bbb H$ is such that $\<\Bbb G+a\>=\Bbb H$;
  \bul $h:\Si\to\Bbb G$ is H\"older continuous and $\sigma$-{\it aperiodic} in the sense that if  $\g\in\widehat{\Bbb G},\ \l\in\Bbb S^1,\ g:\Si\to\Bbb S^1$ H\"older continuous and $\g\circ f=\l\tfrac{g\circ \sigma}g$, then $\g\equiv 1\ \l=1$ and $g$ is  constant.
  \

  It follows that $\dim(\Bbb G)=\dim(\Bbb H)=:D$ where for $A\subset\Bbb R^d$, $\dim(A)$ denotes the dimension of the closed linear subspace spanned by $A$.
  \

  Any closed subgroup $\Bbb G\le\Bbb R^{D}$ with $\dim(\Bbb G)=D$ is conjugate by linear map  to a group of form  $\Bbb Z^k\x\Bbb R^{D-k}$ where
  $0\le k\le D:=\dim(\Bbb G)$ and $\Bbb Z^0,\ \Bbb R^0:=\{0\}$.
\

 Now suppose  that $f:\Si\to\Bbb G=\Bbb Z^k\x\Bbb R^{D-k}$ is H\"older continuous and $\sigma$-aperiodic
 and consider the {\it $(\Si,f)$-random walk} $(\Si\x\Bbb G,\B(\Si\x\Bbb G),\widetilde{m},\s_f)$ where
 $\s_f:\Si\x\Bbb G\to\Si\x\Bbb G$  is defined by
 $$\s_f(x,y):=(\s(x),y+f(x))\ \& \\ d\widetilde{m}(x,y):=d\mu(x)dy$$
   where $\mu$ is the $\s$-invariant Parry measure (with maximal entropy) and $dy$ is Haar measure on $\Bbb G$. \

   As shown in [G], by the aperiodicity of $f$,
$(\Si\x\Bbb G,\s_f,\widetilde{m})$ is exact.
 \

\subsection*{Random walk adic transformation over $(\Si,f,\tau)$}
 The {\it random walk adic transformation over} $(\Si,f,\tau)$ is that skew product over $\tau$ which parametrizes the tail $\mathfrak T(\s_f)$ of the $(\Si,f)$-random walk.

 \

 To identify this:
 \begin{align*}\mathfrak T&(\s_f):=\bigcup_{n\ge 0}\{((x,y),(x',y'))\in(\Si\x\Bbb G)^2:\ \s_f^n(x',y')=\s_f^n(x,y)\}\\ &=
 \bigcup_{n\ge 0}\{((x,y),(x',y'))\in(\Si\x\Bbb G)^2:\ \s^n(x')=\s^n(x)\ \&\ y'+f_n(x')=y+f_n(x)\}\\ &=\{((x,y),(x',y'))\in(\Si\x\Bbb G)^2:\ (x,x')\in\frak T(\s)\ \&\ y'=y+\psi(x,x')\}\end{align*}
 where
 $$\psi(x,y):=\sum_{k=0}^\infty (f(\sigma^ky)-f(\sigma^kx)).$$
 Thus
 \begin{align*}\mathfrak T(\s_f)\cap(\Si'\x\Bbb G)^2&=\{((x,y),(x',y'))\in(\Si\x\Bbb G)^2:\ \s_f^n(x',y')=\s_f^n(x,y)\}\\ &=
 \bigcup_{n\ge 0}\{((x,y),T^n(x,y)):\ (x,y)\in(\Si\x\Bbb G)^2\ n\in\Bbb Z\}\end{align*}
 where
 $T:\Si'\x\Bbb G\to\Si'\x\Bbb G$ is   given by
 \begin{align*}\tag{\symqueen}&T(x,y)=\tau_{\phi}(x,y)=(\tau(x),y+\phi(x))\\ &
 \text{ and}\\ &
 \phi(x):=\psi(x,\tau(x))=\sum_{k=0}^\infty (f(\sigma^k\tau x)-f(\sigma^kx)).\end{align*}
 We consider $T$ with the invariant measure
 $$dm(x,y):=d\nu(x)dy$$
 where $\nu\in\mathcal P(\Si)$ is the $\tau$--invariant Markov measure (equivalent to the Parry measure $\mu$) and $dy$ is Haar measure on $\Bbb G$.

 As mentioned above, it was shown in [G]  that
$(\Si\x\Bbb G,\s_f,\widetilde{m})$ is exact, whence $(\Si\x\Bbb G,T,m)$ is ergodic.
\

The {\it degree} of the random walk adic transformation $(\Si\x\Bbb G,T,m)$ is
$$\text{\tt deg}\,(T):=\dim\,(\text{span}\,\Bbb G).$$

In this paper we ignore the other invariant measures for $T$ (which are considered in [ANSS]).

\section*{\S2 Uniform Convergence}
\

\proclaim{Uniform Convergence Lemma}
\

 \ \ \ Let $\Si$ be a mixing TMS and let
 $\tau:\Sigma'\to\Sigma'$ be the associated stationary adic transformation with $\tau$--invariant Borel probability measure $\nu_0\in\mathcal P(\Si)$, then
 \begin{align*}\frac1n\sum_{k=0}^{n-1}F\circ \tau^k\xrightarrow[n\to\infty]{}\int_{\Sigma'} Fd\nu_0\ \ \ \text{uniformly on $\Sigma'$}\ \ \ \forall\ \ F\in C(\Si).
 \end{align*}
 \endproclaim
Although $\tau$ is a uniquely ergodic homeomorphism on $\Si'$, this space  is not compact.
 \

 The main part of the proof is to provide a suitable
continuous transformation of a related  compact space.
This latter construction is made for  any adic transformation.

\

\proclaim{Compact Representation Lemma}
\

Let $\Si$ be a Bratteli diagram  and let $\tau$ be the associated adic transformation.
     There are: \sbul a compact metric space $(\widehat\Sigma,\widehat{d})$;\sbul a  continuous injection  $\pi:\Si\setminus\Si_{\max}\to\widehat\Sigma$ and \sbul continuous surjections $\varpi:\widehat\Sigma\to\Sigma,\ \widehat\tau:\widehat\Sigma\to\widehat\Sigma$ so that
 $$\varpi\circ\pi=\mathrm{Id}|_{\Si\setminus\Si_{\max}},\ \pi\circ\tau=\widehat\tau\circ\pi
\ \&\ \varpi\circ\widehat\tau=\tau\circ\varpi.$$
\endproclaim
\begin{proof}

 For $\om\in\Sigma_{\max}$, set
$$A(\om):=\{\a\in\Sigma:\ \exists\textrm{ sequence } x^{(n)}\in\Sigma\setminus\Si_{\max}\ \textrm{ s.t. }x^{(n)}\to\om\textrm{ and } \tau(x^{(n)})\to\a\},$$
then $A(\om)\subset\Sigma_{\min}$.
\

Let
$$\Si_0:=\Si\setminus\Si_{\max}\ \ \&\ \ \Si_1:=\bigcup_{\om\in\Sigma_{\max}}\{\om\}\x A(\om).$$  Define the metric space $(\widehat\Si,\widehat{d})$ by
$$\widehat\Sigma:=\Sigma_0\uplus\Si_1;\ \ \ \widehat{d}|_{\Si_0\x\Si_0}\equiv d\ \&\ \text{for}\ (\om,\a)\in\Si_1:$$
$$\widehat{d}((\om,\a),z)=
\begin{cases} d(\om,z)+d(\a,\tau z) & z\in\Si_0;\\
d(\om,\om')+d(\a,\a') & z=(\om',\a')\in\Si_1.\end{cases}$$
To see that this is a  compact metric space, let $(z_n)_{n\ge 1}$ be a sequence in $\widehat{\Si}$,
then:
\begin{itemize}
\item If $\exists\ n_k\to\infty,\ z_{n_k}=(\om_k,\a_k)\in\Si_1$, then (possibly passing to a subsequence) we may assume that $(\om_k,\a_k)\to (\om,\a)\in\Si_{\max}\x\Si_{\min}$. To see that $(\om,\a)\in\Si_1$, for each $k\ge 1$, choose $x_k\in\Si_0$ so that $d(x_k,\om_k)+d(\tau(x_k),\a_k)\to 0$. It follows that
$(\om,\a)\in\Si_1$ and $(\om_k,\a_k)\to(\om,\a)$ in $\widehat{\Si}$.
\item Otherwise $\exists\ n_k\to\infty,\ z_{n_k}\in\Si_0,\ z_{n_k}\to r\in\Si$ and
\begin{itemize}
\item  if $r\in\Si_0$, then $z_{n_k}\to r$ in $\widehat{\Si}$;
\item if $r\in\Si_1$, then $\exists\ m_\ell=n_{k_\ell}\to\infty$ so that $\tau(z_{n_k})\to s\in\Si$; whence $(r,s)\in\Si_1$ and $z_{m_\ell}\to(r,s)$ in $\widehat{\Si}$.
\end{itemize}
\end{itemize}

Let $\pi:\Si_0\to \widehat{\Sigma}$ be the identity map. The following map is a continuous left inverse:
 $\varpi:\widehat\Sigma\to\Sigma$ defined by the identity map on $\Si_0$ and by $(\omega,\a)\mapsto
 \omega$ on $\Si_1$.

Now define $\widehat\tau:\widehat\Sigma\to\widehat\Sigma$ by $\widehat\tau(x):=\tau(x)$ for $x\in\Sigma_0$ and $\widehat\tau(\om,\a)=\a$ for $(\om,\a)\in\Sigma_1$, then $\widehat\tau$ is continuous and $\pi\circ\tau=\widehat\tau\circ\pi$.
 \

 To see that $\widehat\tau$ is onto, it suffices to show that $\widehat\tau(\widehat\Si\supset\Si_{\min}$. To this end,  fix $\a\in\Si_{\min}$, then $\exists\ x_n\in\Si\setminus\Si_{\min},\ x_n\to\a$. Without loss of generality, $x_n\notin\Si_{\max}$ and so $x_n=\tau(y_n)$ for some $y_n\in\Si$ where $y_n\to\om\in\Si$. It follows that $\om\in\Si_{\max}$ (else $\a=\tau(\om)\notin\Si_{\min}$) whence $(\om,\a)\in\widehat{\Si}$ and $\varpi$ is onto.
\end{proof}
\demo{Proof of the Uniform Convergence Lemma}
\

Since $\nu_0\in\mathcal P(\Si')$, it lifts to a $\widehat{\tau}$--invariant measure $\nu_1\in\mathcal P(\widehat{\Sigma})$: $\nu_1\circ\pi^{-1}=\nu_0$. We claim that
 $\widehat\tau$ is uniquely ergodic on $\widehat\Sigma$ with invariant measure $\nu_1$. Else $\exists\ \nu_1\ne \nu_2\in\mathcal P(\widehat{\Sigma})$ with $\nu_2\circ\widehat{\tau}^{-1}=\nu_2$. This entails $\nu_2\circ\varpi^{-1}=\nu_0$ whence $\nu_2=\nu_1$.
 \

 It follows that
$$\Delta_n(F):=\sup_{\widehat\Sigma}\bigg|\frac1n\sum_{k=0}^{n-1}F\circ\widehat\tau^k-\int_{\widehat\Sigma}Fd\nu_1\bigg|
\xrightarrow[n\to\infty]{} 0\ \forall\ F\in C(\widehat\Sigma).$$
 If $f\in C(\Sigma)$, then $f\circ\varpi\in C(\widehat\Si)$ and for every $\om\in\Sigma'$ and $k\ge 1$,  we have $\tau^k\om\in\Sigma'$ whence
$f(\tau^k\om)=f\circ\varpi(\widehat\tau^k\pi\om).$ Thus
$$\sup_{\om\in\Si'}\big|\tfrac1n\sum_{k=0}^{n-1}f(\tau^k\om)-\int_{\widehat\Sigma}fd\nu_1\big|
\le\Delta_n(f\circ\phi)\xrightarrow[n\to\infty]{}0.\ \ \ \Checkedbox$$

\section*{\S3 Limit theory for the shift}
\

Let $(\Si,\mathcal A,\mu,\sigma)$ be a mixing TMS with $\mu$ a $\s$-invariant, Markov measure
 let $\widehat\s=\widehat\s_\mu$ be its transfer operator, let $\Bbb G\subset\Bbb R^d$ be a closed subgroup of dimension $D$; and let $f=(f^{(1)},\dots,f^{(D)}):\Si\to \Bbb G$ be  H\"older
continuous  and aperiodic. Let   $\o f:=f-\Bbb
E( f)$. Note that $\o f$ may have values outside $\Bbb G$.
\

We'll need the following results for the sequel. The results are not new although we did not find references for their statements. The proofs are standard and will only be sketched.
\proclaim{3.1\ Asymptotic variance theorem}
\par  $\exists$ a $D\x D$, symmetric, positive definite matrix $\G=\G_f$ so that
$$\text{\tt Rank}\,\G_f=D\ \&\ \frac1n\Bbb E(\o f^{(i)}_n\o f^{(j)}_n)\underset{n\to\infty}\lra\ \G_{i,j}\ \ \forall\ 1\le i,j\le D$$
where $\o f^{(i)}_n:=\sum_{k=0}^{n-1}\o f^{(i)}\circ \sigma^k$.\endproclaim

Since $\Gamma$ is symmetric and positive definite, it can be put in the form
$$
\Gamma=(U^t M)(U^t M)^t
$$
with $U$ unitary and $M>0$ diagonal.

\proclaim{3.2 Central limit theorem}
\begin{equation*}\tag{\tt CLT}\widehat{\sigma}^n1_{[ \frac{ \o f_n}{\sqrt n}\in I]}\underset{n\to\infty}\lra\ \frac1{\sqrt{(2\pi)^D\det\G}}\int_I\exp[{-\tfrac12 u^t\G^{-1} u}]du\end{equation*} whenever $I\subset\Bbb R^D$ is Riemann integrable (i.e. with Riemann integrable indicator function).\endproclaim

\

\proclaim{3.3  Local limit theorem\ \ {\rm(mixed lattice-nonlattice case)} }

\

Let  $0\in I\subset\Bbb G$ be  Riemann integrable; then
\begin{equation*}\tag{\tt LLT}(\sqrt n)^D\widehat{\sigma}^n1_{[ f_n\in n\Bbb E(f)+\sqrt nx_n+I]}
\xrightarrow[n\to\infty,\ x_n\to u]{}\frac{m_\Bbb G(I)}{\sqrt{(2\pi)^D\det\G}}\cdot
\exp[{-\tfrac12 u^t\G^{-1} u}].\end{equation*}\endproclaim
\demo{Proof sketches}
\

Suppose that $f$ is $\th$-H\"older continuous and  let $\rho=\rho_\th$ be the metric on $\Si$ defined by
$$\rho(x,y):=\th^{\inf\,\{n\ge 1:\ x_n\ne y_n\}}.$$ Let $\Bbb L=\Bbb L_\th$ be the Banach space of $\rho$-Lipschitz continuous (equivalently, $\th$-H\"older continuous functions) on $\Si$ equipped with the norm
$$\|F\|_\Bbb L:=\sup_{x\in\Si}|F(x)|+\sup_{x,\ y\in\Si}\frac{|F(x)-F(x)|}{\rho(x,y)}.$$
As shown in [D-F], [GH], for some $N\ge 1,\ \widehat\s^N:\Bbb L\to\Bbb L$  satisfies the {\it Doeblin-Fortet inequality}, namely $\exists\ r\in (0,1)\ \&\ H>0$ so that
\begin{equation*}\tag{\tt DF}\|\widehat{\sigma}^NF\|_\Bbb L\le r\|F\|_\Bbb L+H\sup_{x\in\Si}|F(x)|\ \ \forall\ F\in\Bbb L,\end{equation*}
 whence ([D-F], [IT-M]) $\widehat\s:\Bbb L\to\Bbb L$ is
quasi-compact in that $\exists\ \th\in (0,1)\ \&\ M>1$ so that
\begin{equation*}\tag{\tt QC}\|\widehat\s^nF-\int_{\Si}Fd\mu\|_\Bbb L\le M\th^n\|F\|_\Bbb L\ \ \forall\ F\in\Bbb L.\end{equation*}
\demo{Proof of the asymptotic variance theorem}
\

The  absolute convergence of the series
$\sum_{n\ge 1}|\int_{\Si}f^{(i)}\cdot f^{(j)}\circ\s^nd\mu|$ for $1\le i,j\le D$ follows from ({\tt QC}) and  the convergence follows from this. The non-singularity of the limit matrix follows from the aperiodicity of each $t\cdot f:\Si:\to\Bbb R\ \ \ (t\in\Bbb R^D)$ via Leonov's theorem (as in [RE]).\ \ \Checkedbox

\

\demo{Proof of the central and local limit theorems}
\

Both of these results are established   using Nagaev's perturbation method as in [HH], [PP], [RE] (aka characteristic function operators [AD]).
\

For $t\in\Bbb C^D$ define $P_t:L^1(\Si)\to L^1(\Si)$ by
$P_t(F):=\widehat\s(e^{it\cdot f}F),$ then
$P_t(F):\Bbb L\to\Bbb L$ and the map $t\mapsto P_t$ is analytic $\Bbb C^D\to\hom(\Bbb L,\Bbb L)$ with
$$(\frac{\bdy^r\ \ \ }{\bdy t_{k_1}\cdots\bdy t_{k_r}}P_t)(F):=i^r\widehat\s(\prod_{j=1}^rf^{(k_j)}e^{it\cdot f}F);$$
where $\hom(\Bbb L,\Bbb L)$ is equipped with the uniform topology. We have (see [N], [PP] $\&$/or [RE])
\bul $\|P_t\|\le 1\ \ \forall\ t\in\widehat{\Bbb G}=\Bbb T^k\x\Bbb R^{D-k}$ with equality iff $t=0$;
\bul $P_t$ satisfies ({\tt DF}) for $|t|$ small; whence
\proclaim{ Nagaev's Theorem  \ \ [N]}\ \ \  There are constants $\e>0,\ K>0$ and
$\th\in (0,1)$; and analytic functions $\l:(-\e,\e)\to B_{\Bbb
C}(0,1),\ N:(-\e,\e)\to \hom(L,L)$ such that
$$\|P_t^nh-\l(t)^n N(t)h\|_L\le K\th^n\|h\|_L
 \ \ \ \ \forall\ |t|<\e,\ n\ge 1,\ h\in L$$
where $\forall |t|<\epsilon$, $N(t)$ is a projection onto a
one-dimensional subspace, $\l(0)=1\ \&\  N(0)F:=\int_\Si Fd\mu$.\endproclaim
\

The expansion of $\l$ is  obtained by considering $t\cdot f:\Si\to\Bbb R \ \ (t\in\Bbb R^D)$ as in  [GH]. It is given by
 $$\l(t)=1+it\cdot\Bbb E(f)-\frac{<\G t,t>}2+o(|t|^2)\ \ \text{as}\ t\to 0.$$

 The central  limit theorem follows from this in the standard manner (see [GH], [RE]); and the local limit theorem  follows with a proof as in [S] (see [AD]).

\section*{\S4 Proof of the theorem }
\

For $n\ge 1$, set $\ell_n:=[\log_\l n]$. Let $0\in I\subset\Bbb G$ be Riemann integrable with $0<|I|<\infty$ of form $I=\{0^{(k)}\}\x J$ where $0^{(k)}\in\Bbb Z^k,\ 0^{(k)}_j=0\ (1\le j\le k)$ and  $0\in J\subset\Bbb R^{D-k}$ is Riemann integrable with $0<|J|<\infty$.
\

To achieve our goal, we'll establish:

\begin{align*}\tag{\Football}&\forall\ R>0,\ \forall\ x\in\Si',\\ &
\frac{\ell_n^{D/2}}{n}\cdot S_n(1_{\Si'\x I})(x,0)1_{B(R)}(\tfrac{\o f_{\ell_n}(x)}{\sqrt{\ell_n}})\ \ \approx\ \ \\ &\ \ \ \ \ \ \ \ \ \ \frac{|I|}{\sqrt{(2\pi)^D\det\G}}\cdot\exp[{-\tfrac{\|M^{-1}U\o f_{\ell_n}(x)\|^2}{2\ell_n}}]1_{B(R)}(\tfrac{\o f_{\ell_n}(x)}{\sqrt{\ell_n}})\end{align*}
where $\o f:=f-\Bbb E(f)$ and $B(R):=\{z\in\Bbb R^D:\ \|z\|<R\}$ and
$a_n\approx b_n$ means $a_n-b_n\underset{n\to\infty}\lra\ 0$.

We'll show first that {\rm (\Football)}  holds, and then we'll prove that (\Football)  $\Lra$ the theorem.

\subsection*{Overview of the proof of (\Football)}
 \

The proof uses a process of  {\tt block splitting} where in order to estimate
$$S_n^{(T)}(1_{\Om\x I})(x,0)=\sum_{j=0}^{n-1}1_{\Om\x I}(\tau^kx,\phi_k(x))$$ we split the $\tau$-orbit block
  $\{\tau^kx:\ 0\le k\le n-1\}$ into simpler blocks on which it is easy to apply the results of \S3.

  \

  This is done as follows.
  \

 For $x\in\Om,\ N\ge 1$
$$\s^{-N}\{\s^Nx\}=\{\tau^kx_{\text{\tt\tiny min}}: 0\le k\le \#\s^{-N}\{\s^Nx\}-1\}$$
where $x_{\text{\tt\tiny min}}:=\min\,\s^{-N}\{\s^Nx\}$ with respect to  the reverse lexicographic order and
$$\sum_{j=0}^{\#\s^{-N}\{\s^Nx\}-1}1_{\Om\x I}\circ T^k(x_{\text{\tt\tiny min}},0)=\#\{y\in\s^{-N}\{\s^Nx\}:\ f_N(y)\in f_N(x_{\text{\tt\tiny min}})+I\}.$$
Quantities appearing, such as $$\#\{y\in\s^{-N}\{x\}:\ \ f_N(y)\in f_N(z)+I\}$$ where $I\subset\Bbb G$ is Riemann integrable, are estimated using ({\tt LLT}) as in lemma 4.1 (below).
\

The arbitrary blocks are estimated  from the simple ones of suitably smaller size. This is calibrated using
$$\#\s^{-N}\{\s^Nx\}=\sum_{0\le s\le d-1}A^N_{s,x_{N+1}}\ \sim\ c(x_{N+1})\l^N$$
where $A$ is  transition matrix of $\Si$ and $\l=e^{h_{\text{\tiny top}}(\Si)}$  is its leading  eigenvalue.

\

 Fix $M\ge 1$ large. For each $n\ge 1$ large, let $N=N_n$ be such that
$M\l^N=\l^{\pm 1} n$, then
 $$\{\tau^jx:\ 0\le j<n\}=\bigcupdot_{k=0}^{M-1}\s^{-N}\{\tau^k\s^{N}x\}$$ up to relatively small  edge effectss (estimated in the proof below using lemma 4.2) and

\begin{align*}S_n^{(T)}&(1_{\Om\x I})(x,0)=\\ &\sum_{k=0}^{M-1}\#\{y\in\s^{-N}\{\tau^k\s^Nx\}:\ f_N(y)=f_N(\min\,\s^{-N}\{\tau^k\s^Nx\})\}\end{align*} up to relatively  small errors (estimated in lemma 4.3 below).
\

\subsection*{Proof of (\Football)}
\

For $x\in\Si,\ t^{(n)}\in\Bbb G,\  \ \sup_n\tfrac{\|t^{(n)}\|}{\sqrt n}<\infty$, set $$N_n(x):=\#\{z\in\s^{-n}\{\s^nx\}:\ \ \o f_n(z)\in t^{(n)}+I\}.$$
\proclaim{Lemma 4.1}
 \begin{align*}N_n(x)\ \ \sim \ \ \frac{\l^nh(\sigma^nx)|I|}{\sqrt{(2\pi n)^{D}\det\G}}\exp[{-\frac{\|M^{-1}Ut^{(n)}\|^2}{2n}}]\end{align*} \ \ {\it uniformly on $\Si$ where $h=\frac{d\mu}{d\,m}$;
 $m$ and  $\mu$ being the $\tau$- and  $\sigma$-invariant probabilities (respectively).}\endproclaim \demo{Proof}

 \

 \ \ \  Let $\widehat \sigma_m,\ \widehat \sigma_{\mu}$ be the transfer operators of $\s$ with respect to $m\ \&\ \mu$ respectively, then $\widehat \sigma_mf=h\widehat \sigma_{\mu}(\tfrac{f}h)$ and

\begin{align*}N_n(x)&=\l^n\widehat \sigma_m^n1_{[\o f_n\in t^{(n)}+I]}(\sigma^nx)\\ &=\l^nh(\sigma^nx)\widehat \sigma_\mu^n(\tfrac{1_{[\o f_n\in t^{(n)}+I]}}h)(\sigma^nx)\\  &\sim \frac{\l^nh(\sigma^nx)|I|}{{\sqrt{(2\pi n)^{D}\det\G}}}\exp[{-\frac1{2n} \|M^{-1}Ut^{(n)}\|^2}]\end{align*}
\ \  uniformly on $\Si$ by ({\tt LLT}).\ \ \Checkedbox

\subsection*{Block splitting}
\

 For $n\ge 1$, let $\Si_{n,s}:=\{(x_1,\dots,x_n):\ x\in\Si,\ x_{n+1}=s\},\ J_n(s):=\#\Si_{n,s}$, then
$$J_n(s)=\sum_{u\in \mathcal S}A^n_{u,s}\underset{n\to\infty}\sim\ c(s)\l^n$$  uniformly in $s\in\mathcal S$ where $\l=e^{h_{\text{\tiny top.}}(\Si,\sigma)}$ and $c:\mathcal S\to\Bbb R_+$.
\

It will be convenient also to set $\widehat J_n(z):=\#\,\s^{-n}\{z\}$. Here
$$\widehat J_n(x)=\#\Si_{n,x_{1}}=J_n(x_{1})$$ and
$$\widehat J_n(x)\sim\ \frak c(x)\l^n$$  uniformly in $x\in\Si$ where $\frak c:\Si\to\Bbb R_+,\ \frak c(x):=c(x_1)$.

\

For $n\ge 1$ fixed, we call a point $x\in\Si$
\bul $n$-{\it minimal} if $x=\min\,\sigma^{-n}\{\sigma^nx\}=\min\,\{y\in \Si:\ y_{n+1}=x_{n+1}\}\ \&$
\bul $n$-{\it maximal} if $x=\max\,\sigma^{-n}\{\sigma^nx\}=\max\,\{y\in \Si:\ y_{n+1}=x_{n+1}\}$;
\

Now define $$K_n:\Si\to\Bbb N\ \&\ \tau_n:\Si\to\Si\ \ \text{\rm by}$$
$$K_n(x):=\min\,\{k\ge 1:\ \tau^kx\ \text{ is $n$-maximal}\}\ \&\ \tau_n(x):=\tau^{K_n(x)+1},$$
then:
\bul $\sigma^n\tau_n(x)=\tau(\sigma^nx)$;\bul $\tau_n(x)$ is $n$-minimal and
\bul $K_n(x)\le \widehat J_n(\s^nx)=\#\sigma^{-n}\{\sigma^nx\}$ with equality if $x$ is $n$-minimal.

\

It follows that for $j\ge 1$,
$$\s^n\tau_n^{j}(x)=\tau^j(\s^nx)$$
and
$$K_n(\tau_n^{j}(x))=\widehat J_n(\tau^j\s^nx).$$

\

For $n\ge 1$ fixed and $r\ge 1$, set
$$K_n^{(r)}(x):=\sum_{j=0}^{r-1}K_n(\tau_n^j(x))=K_n(x)+\sum_{j=1}^{r-1}\widehat J_n(\tau^j(\s^nx)).$$
\proclaim{Lemma 4.2}\ \ \ \  $\exists\ \eta_n,\ \theta_r\downarrow\ 0$ so that
$$K_n^{(r)}(x)=e^{\pm(\eta_n+\theta_r)}r\l^nE(\frak c)\ \ \forall\ n,r\ge 1,\ x\in\Si'.$$\endproclaim
\demo{Proof}\ \ \   By  the uniform convergence lemma  $\exists\ \th_r\downarrow 0\ \st$
 $$\sum_{j=1}^{r-1}\frak c(\tau^j(\s^nx))=e^{\pm\theta_r}rE(\frak c)\ \forall\ x\in\Si',\ n,\ r\ge 1.$$

Suppose that $J_n(s)=e^{\pm\eta_n}\l^n\frak c(s)$ where $\eta_n\downarrow 0$, then for $x\in\Si'$,
\begin{align*}K_n^{(r)}(x)&=K_n(x)+\sum_{j=1}^{r-1}\widehat J_n(\tau^j(\s^nx))\\ &=K_n(x)+e^{\pm\eta_n}\l^n\sum_{j=1}^{r-1}\frak c(\tau^j(\s^nx))\\ &=K_n(x)+e^{\pm\eta_n}e^{\pm\theta_r}r\l^nE(\frak c).\end{align*}
Since $K_n(x)\le \widehat J_n(\s^nx)=O(\l^n)$, the lemma follows.\ \ \ \Checkedbox

\proclaim{Lemma 4.3} For $r\in\Bbb N$ fixed, $x\in\Si'$ and $R>0$, as $n\to\infty$:
\begin{align*}\tag{1}S_{K_n^{(r)}(x)}&(1_{\Si\x I})(x,0)1_{B(R)}(\tfrac{\o f_{n}(x)}{\sqrt{n}})\ \gtrsim\ \\ & \tfrac{h_{r-1}(\sigma^nx)|I|\l^n}{n^{D/2}}\exp[{-\tfrac{\|M^{-1}U\o f_n(x)\|^2}{2n}}]1_{B(R)}(\tfrac{\o f_{n}(x)}{\sqrt{n}})\end{align*} and
\begin{align*}\tag{2}S_{K_n^{(r)}(x)}&(1_{\Si\x I})(x,0)1_{B(R)}(\tfrac{\o f_{n}(x)}{\sqrt{n}})\ \lesssim\ \\ & \tfrac{h_{r+1}(\sigma^nx)|I|\l^n}{n^{D/2}}\exp[{-\tfrac{\|M^{-1}U\o f_n(x)\|^2}{2n}}]1_{B(R)}(\tfrac{\o f_{n}(x)}{\sqrt{n}})\end{align*}
where $h_r(z):=\sum_{j=0}^{r-1}h(\tau^jz)$.\endproclaim

\

Here, for $A_n,\ B_n>0,\ A_n\ \gtrsim\ B_n$ means $\varliminf_{n\to\infty}\frac{A_n}{B_n}\ge 1$ and
$\lesssim\ B_n$ means $\varlimsup_{n\to\infty}\frac{A_n}{B_n}\le 1$

\demo{Proof}
\

 Writing $K_n^{(0)}\equiv 0$, we have
\begin{align*}\tag{\Rightscissors} &S_{K_n^{(r)}(x)}(1_{\Si\x I})(x,0)\\ &=\sum_{j=0}^{r-1}\left(S_{K_n^{(j+1)}(x)}(1_{\Si\x I})(x,0)-S_{K_n^{(j)}(x)}(1_{\Si\x I})(x,0)\right)\\ &=S_{K_n(x)}(1_{\Si\x I})(x,0)+\\ & \ \ \ \ \ +\sum_{j=1}^{r-1}\left(S_{K_n^{(j+1)}(x)}(1_{\Si\x I})(x,0)-S_{K_n^{(j)}(x)}(1_{\Si\x I})(x,0)\right).\end{align*}

For fixed $j\ge 1$,
\begin{align*}& S_{K_n^{(j+1)}(x)}(1_{\Si\x I})(x,0)-S_{K_n^{(j)}(x)}(1_{\Si\x I})(x,0)\\ &=
\sum_{k=K_n^{(j)}(x)}^{K_n^{(j+1)}(x)-1}1_{\Si\x I}(\tau^kx, \phi_k(x))\\ &=
\sum_{\ell=0}^{K_n(\tau_n^j(x))-1}1_{\Si\x I}(\tau^\ell(\tau^{K_n^{(j)}(x)}x),\phi_{K_n^{(j)}(x)+\ell}(x))\\ &=\sum_{\ell=0}^{\widehat J_n(\tau^j\s^n x)-1}1_{\Si\x I}(\tau^\ell(\tau_n^j(x)),\phi_{K_n^{(j)}(x)+\ell}(x))
\end{align*}

Now,
$$\{\tau^\ell(\tau_n^j(x):\ 0\le \ell\le \widehat J_n(\tau^j\s^n x)-1\}= \sigma^{-n}\{\tau^j\s^n x\},$$
$$\phi_{K_n^{(j)}(x)+\ell}(x) =\psi(z,x)=\sum_{k=0}^\infty (f(\sigma^kx)-f(\sigma^kz))\ \ \ (0\le \ell\le \widehat J_n(\tau^j\s^n x)-1);$$
and so \begin{align*} S_{K_n^{(j+1)}(x)}(1_{\Si\x I})(x,0)&-S_{K_n^{(j)}(x)}(1_{\Si\x I})(x,0)\\ &=
\#\{z\in \sigma^{-n}\{\sigma^n\tau_n^j(x)\}:\ \psi(z,x)\in I\}.\end{align*}
For $z=\tau^\ell(\tau_n^j(x)\in \sigma^{-n}\{\tau^j\s^n x\}$ we have
\begin{align*}\phi_{K_n^{(j)}(x)+\ell}(x)) &=\psi(z,x)\\ &=\sum_{k=0}^\infty (f(\sigma^kx)-f(\sigma^kz))\\ &=f_t(x)-f_t(z)\end{align*}
where $t=t(x,z):=\min\,\{N\ge 1:\ \sigma^Nx=\sigma^Nz\}$.
\

Here, $\s^nz=\tau^j\s^nx$ and
$$t(x,z)\le n+t(\tau^j\sigma^nx,\sigma^nx).$$

 Thus

\begin{align*}\psi(z,x) =f_n(z)-f_n(x)+\kappa_{n,j}(x)\end{align*} where
$$|\kappa_{n,j}(z)|\le 2\sup|f|(N-n)\le 2\sup|f|t(\tau^j\sigma^nx,\sigma^nx).
$$
We claim that for a.e. $x\in\Si$,
\begin{align*}\tag{\dsmathematical}\max_{1\le j\le r}\,t(\tau^j\sigma^nx,\sigma^nx)=O(\log n)\ \text{  as }\ n\to\infty.\end{align*}\demo{Proof of (\dsmathematical)}
\

For $M>0$ set $A_n(M):=\{x\in\Si:\ t(\tau\sigma^nx,\sigma^nx)>M\log n\}$, then
\begin{align*}m(\{x\in\Si:\ \max_{1\le j\le r}\,t(\tau^j\sigma^nx,\sigma^nx)>M\log n\})&\le m(\bigcup_{0\le j\le r-1}\tau^{-j}A_n(M))\\ &\le rm(A_n(M)).\end{align*}
Now $t(\tau\sigma^nx,\sigma^nx)>M\log n$ iff $\exists\ z\in\Si_{\max}$ so that
$$x_{n+j}=z_j\ \forall\ 1\le j\le M\log n.$$
Thus
\begin{align*}m(A_n(M))&\le \sup\,h\mu(A_n(M))\\ &= \sup\,h\sum_{z\in\Si_{\max}}\mu([z_1,\dots,z_{\lfl M\log n\rfl}])\\ &=O(\l^{-M\log n})\end{align*}
and $\sum_{n\ge 1}m(A_n(M))<\infty$ for $M>\frac1{\log\l}$. The claim (\dsmathematical) now follows from the Borel-Cantelli lemma.\ \ \Checkedbox

\

In view of (\dsmathematical), we have by lemma 4.1 that for a.e. $x$:
\begin{align*}& 1_{B(R)}(\tfrac{\o f_{n}(x)}{\sqrt{n}})(S_{K_n^{(j+1)}(x)}(1_{\Si\x I})(x,0)-S_{K_n^{(j)}(x)}(1_{\Si\x I})(x,0))\\ &=
1_{B(R)}(\tfrac{\o f_{n}(x)}{\sqrt{n}})\#\{z\in \sigma^{-n}\{\sigma^n\tau_n^j(x)\}:\ f_n(z)\in f_n(x)-\kappa_{n,j}(x)+I\}\\ &\sim
|I|1_{B(R)}(\tfrac{\o f_{n}(x)}{\sqrt{n}})\tfrac{\l^nh(\sigma^nx)}{n^{D/2}}\exp[{-\tfrac{\|M^{-1}U(\o f_{n}(x)-\kappa_{n,j}(x))\|^2}{2n}}]\\ &\sim |I|1_{B(R)}(\tfrac{\o f_{n}(x)}{\sqrt{n}})\tfrac{\l^nh(\sigma^nx)}{n^{D/2}}\exp[{-\tfrac{\|M^{-1}U\o f_{n}(x)\|^2}{2n}}].\end{align*}

\

The lemma follows from this.\ \ \Checkedbox

\demo{Completion of  the proof of (\Football)}

 Given $0<\e<\tfrac13$,

 \bul use the uniform convergence lemma  to fix $r_\e\ \st\ \forall\ y\in\Si',\ r\ge r_\e,\ \ e^{\th_r}<1+\e$  where $\th_r$ is as in lemma 4.2, and
 $$\frak c_r(y),\ \frak c_{r+2}(y)=(1\pm \e)r\Bbb E_m(\frak c)\ \&\ h_r(y),\ h_{r+2}(y)=(1\pm \e)r\Bbb E_m(h).$$
 \bul fix $J>e^{\th_r}\ \forall\ r\ge 1$ and for $n\ge 1$ let
 $$L_n=L_{n,\e}:=\lfl\log_\l\tfrac{n}{2JE(c)r_\e}\rfl$$
 and let $r_n=r_{n,\e}$ be $\st$
 $$K_{L_n}^{(r_n)}(\tau \sigma^nx)\le n<K_{L_n}^{(r_n+1)}(\tau \sigma^nx)<K_{L_n}^{(r_n+2)}(\sigma^nx).$$
 It follows that
 $$S_{K_{L_n}^{(r_n)}(\tau \sigma^nx)}(1_{\Si\x I})(x,0)\le S_{n}(1_{\Si\x I})(x,0)\le S_{K_{L_n}^{(r_n+2)}(\sigma^nx)}(1_{\Si\x I})(x,0).$$
 By lemma 4.2,
 $$n\le K_{L_n}^{(r_n+2)}(\sigma^nx)\le e^{(\eta_{L_n}+\theta_{r_n})}r_n\l^{L_n}E(c)\lesssim e^{\theta_{r_n}}r_n\tfrac{n}{2JE(c)r_\e}E(c)\le \tfrac{nr_n}{2r_\e}$$
 whence for large $n$, $r_n>r_\e$  and
 \begin{align*}S_{n}(1_{\Si\x I})(x,0)1_{B(R)}(\tfrac{\o f_{L_n}(x)}{\sqrt{L_n}})&\ge S_{K_{L_n}^{(r_n)}(\tau \sigma^nx)}(1_{\Si\x I})(x,0)1_{B(R)}(\tfrac{\o f_{L_n}(x)}{\sqrt{L_n}})
 \\ &\gtrsim\ \ \tfrac{h_{r_n}(\tau \sigma^{L_n}x)\l^{L_n}}{{L_n}^{D/2}}\exp[{-\tfrac{\|M^{-1}U\o f_{L_n}(x)\|^2}{2{L_n}}}]1_{B(R)}(\tfrac{\o f_{{L_n}}(x)}{\sqrt{{L_n}}})\\ &\ge
 (1-\e)|I|\tfrac{\l^{L_n}r_n}{{L_n}^{D/2}}\exp[{-\tfrac{\|M^{-1}U\o f_{L_n}(x)\|^2}{2{L_n}}}]1_{B(R)}(\tfrac{\o f_{{L_n}}(x)}{\sqrt{{L_n}}})\end{align*}
 and similarly
 \begin{align*}S_{n}(1_{\Si\x I})(x,0)1_{B(R)}(\tfrac{\o f_{L_n}(x)}{\sqrt{L_n}})&\le S_{K_{L_n}^{(r_n+2)}(\sigma^nx)}(1_{\Si\x I})(x,0)1_{B(R)}(\tfrac{\o f_{L_n}(x)}{\sqrt{L_n}})
 \\ &\lesssim\ \ \tfrac{h_{r_n+2}(\sigma^{L_n}x)\l^{L_n}}{{L_n}^{D/2}}\exp[{-\tfrac{\|M^{-1}U\o f_{L_n}(x)\|^2}{2{L_n}}}]1_{B(R)}(\tfrac{\o f_{{L_n}}(x)}{\sqrt{{L_n}}})\\ &\le
 (1+\e)|I|\tfrac{\l^{L_n}r_n}{{L_n}^{D/2}}\exp[{-\tfrac{\|M^{-1}U\o f_{L_n}(x)\|^2}{2{L_n}}}]1_{B(R)}(\tfrac{\o f_{{L_n}}(x)}{\sqrt{{L_n}}}).\end{align*}

 Now,
 \bul $L_n\sim \ell_n=\log_\l n\ \forall\ \e>0$ and since $r_n>r_\e$,
 \begin{align*}n\ge K_{L_n}^{(r_n)}(\tau \sigma^nx)\ge e^{-(\eta_n+\theta_{r_n})}r_n\l^{L_n}E(c)\gtrsim \tfrac1{1+\e}r_n\l^{L_n}E(c).\end{align*}
  \begin{align*}n<K_{L_n}^{(r_n+1)}(\tau \sigma^nx)\lesssim (1+\e)r_n\l^{L_n}E(c)<\tfrac1{1-\e}r_n\l^{L_n}E(c)\end{align*}
  whence
  $$\tfrac{\l^{L_n}r_n}{{L_n}^{D/2}}=(1\pm\e)^2\tfrac{n}{\ell_n^{D/2}E(c)}.$$
This proves (\Football).\ \ \ \ \Checkedbox

\

  \demo{Proof that (\Football) $\Lra$ the theorem.} \ \  \    Let $\G$ be as in \S3 and write $\G=VV^t$ where $V:=U^tM$ with $U$  unitary and $M>0$ diagonal. \sms Let $\xi=(\xi_1,\dots,\xi_d)$ where $\xi_1,\dots,\xi_d$ are independent, identically distributed, Gaussian random variables with $E(\xi_j)=0,\ E(\xi^2_j)=1$, then
 $Z:=U^tM\xi=V\xi$  is Gaussian with correlation matrix
$$E(Z_iZ_j)=E(\sum_{s,t}V_{i,s}\xi_sV_{j,t}\xi_t)=\sum_{s}V_{i,s}V_{j,s}=\G_{i,j}.$$
 By ({\tt CLT})  $$\tfrac{\o f_n}{\sqrt n}\overset{\mathfrak d}\lra\ U^tM\xi=:Z.$$
 Now suppose that
 (\Football) holds. We'll show that for some $a(n)\propto\tfrac{n}{\ell_n^{D/2}}$ and  for  $g\in C([0,\infty]),\ f\in L^1(m)_+$,
 $$g(\tfrac1{a(n)}\cdot S_n(f))\underset{n\to\infty}\lra\ \Bbb E(g(e^{-\frac12\chi_{D^2}}\cdot m(f)))\ \ \text{weak $*$ in}\ L^\infty(m).$$

 By the asymptotic variance theorem,
 $\Bbb E(\|\o f_n\|_2)=O(\sqrt n)$ and $\forall\ \e>0\ \exists\ R$ so that $m_\Om(A_n(R))>1-\e\ \forall\ n\ge 1$ where $A_n(R):=[\tfrac{\o f_{\ell_n}(x)}{\sqrt{\ell_n}}\in B(R)]$.
 \

 Thus for $n\in\Bbb N\ \&\ R>0$ both large enough and $x\in A_n(R)$ we have
\begin{align*}g(\tfrac{\ell_n^{D/2}}{n}S_n(1_{\Si'\x I})(x,0))= g(\tfrac{|I|}{\sqrt{(2\pi)^D\det\G}}\cdot\exp[{-\tfrac{\|M^{-1}U\o f_n(x)\|^2}{2n}}])\pm\e.\end{align*}
Next, by ({\tt CLT}),
\begin{align*}\int_{\Om\x I}g(\tfrac{|I|}{\sqrt{(2\pi)^D\det\G}}\cdot\exp[{-\tfrac{\|M^{-1}U\o f_n\|^2}{2n}}])dm &
\underset{n\to\infty}\lra \Bbb E(g(\tfrac{|I|}{\sqrt{(2\pi)^D\det\G}}\cdot\exp[{-\tfrac{\|M^{-1}UZ\|^2}{2}}]))\\ &=\Bbb E(g(\tfrac{|I|}{\sqrt{(2\pi)^D\det\G}}\cdot \exp[-\tfrac{\chi_D^2}{2}])).\end{align*}
Thus,  $\exists\ a(n)\propto\tfrac{n}{\ell_n^{D/2}}$,
 $$\int_{\Om\x I}g(\tfrac1{a(n)}\cdot S_n(1_{\Si'\x I}))dm \underset{n\to\infty}\lra\ \Bbb E(g(m_\Bbb G(I)\cdot 2^{D/2}e^{-\frac{1}{2}\chi_{D}^2})).$$
 Using Corollary 3.6.2 of [A1], we obtain that $\forall\ F\in L^1(m)_+,\ g\in C([0,\infty])$,
 $$g(\tfrac1{a(n)}\cdot S_n(F))\underset{n\to\infty}\lra\ \Bbb E(g(m(F)\cdot 2^{D/2}e^{-\frac{1}{2}\chi_{D}^2}))\ \ \text{weak $*$ in}\ L^\infty(m)$$ where $m(F):=\int_{\Si\x\Bbb G}Fdm$.
\ \ \ \ \ \Checkedbox

\section*{\S5 Application to exchangeability}
\

Let $\mathcal S=\{0,1,\dots,d-1\}$ and let $\Si\subset\mathcal S^\Bbb N$ be a mixing TMS.
Define  $F^{\natural}:\Si\to\Bbb Z^{d-1}$ by $F^{\natural}(x)_k:=\d_{x_1,k}$.

As shown in [ADSZ],  $F^{\natural}:\Si\to\Bbb Z^{d-1}$ is $\s$-aperiodic iff $\Si$ is {\it almost onto} in the sense that

$\forall\ b,c\in\mathcal S,\ \exists\  n\ge 1,\ b=a_0,a_1,\dots,a_n=c\in\mathcal S$
  such that  $$\s[a_k]\cap \s[a_{k+1}]\ne\emptyset\ \ (0\le k\le n-1).$$

  Define $\vf:\Si\to\Bbb N$  and $R:\Si\to\Si$ by $$\vf(x):=\min\,\{n\ge 1:\ \tau^n(x)_i=x_{\s(i)}\ \text{\tt some finite permn. $\s$ of $\Bbb N$}\}$$ and $$R(x):=\tau^{\vf(x)}(x).$$

  \

  \proclaim{ Corollary 5.1}
     \

    Suppose that $\Si$ is almost onto, then
    $(\Si,\B(\Si),R,m)$ is an ergodic, probability preserving transformation and $\exists\ b(n)\propto \tfrac{n}{(\log n)^{(d-1)/2}}\ \st$
 $$\frac1{b(n)}\sum_{k=0}^{n-1}\vf\circ R^k\overset{\mathfrak d}\lra\ e^{\frac{1}{2}\chi^2_{d-1}}.$$\endproclaim
 \subsection*{Proof}
\

The random walk adic $$(\Si\x\Bbb Z^{d-1},\B(\Si\x\Bbb Z^{d-1}),m\x m_{\Bbb Z^{d-1}},T)$$ over $(\Si,F^{\natural},\tau)$ is conservative and ergodic. Calculation shows that $T_{\Si\x\{0\}}(x,0)=(Rx,0)$ whence
$(\Si,\B(\Si),R,m)$ is an ergodic, probability preserving transformation.
\

By the theorem, $$\frac1{a(n)}S^{(T)}_n(f)\overset{\mathfrak d}\lra \ e^{-\frac{1}{2}\chi^2_{d-1}}\mu(f)\ \forall\ f\in L^1_+$$ where
 $a(n)\propto\frac{n}{(\log n)^{(d-1)/2}}$ (we absorbed the factor in $a(n)$).
\

In particular,
$$\frac1{a(n)}S^{(T)}_n(1_{\Si\x\{0\}})\overset{\mathfrak d}\lra \ e^{-\frac{1}{2}\chi^2_{d-1}}$$
whence by inversion (proposition 1 in [A2]),
$$\frac1{b(n)}\sum_{k=0}^{n-1}\vf\circ R^k\overset{\mathfrak d}\lra\ e^{\frac{1}{2}\chi^2_{d-1}}$$
where $b(n)=a^{-1}(n)\propto \tfrac{n}{(\log n)^{(d-1)/2}}$.\ \ \ \Checkedbox

\section*{\S6 Chi squared laws for horocycle flows}
Let $M_0$ be a compact, connected, orientable, smooth, Riemannian surface with negative sectional curvature, and let $T^1 M_0$ denote the set of unit tangent vectors to $M_0$.
The {\em geodesic flow} on $T^1 M_0$ is the flow which moves a vector $\un{v}\in T^1 M$ at unit speed along its geodesic.

Margulis \cite{Mrg} and Marcus \cite{Mrc} constructed a continuous flow $h^t: T^1 M\to T^1 M$ such that
\begin{enumerate}
\item[(a)] The $h$--orbit of $\vec{v}\in T^1 M_0$ equals $$
W^{ss}(\vec{v}):=\{\vec{u}\ |\ \mathrm{dist}(g^s(\vec{v}),g^s(\vec{u}))\xrightarrow[s\to\infty]{}0\}$$
\item[(b)] $\exists\mu$ s.t. $g^{-s}\circ h^t\circ g^s=h^{\mu^s t}$
\end{enumerate}
In the special case when  $M_0$ is a hyperbolic surface, $h$ is the {\em stable horocycle flow}. Properties (a) and (b) should be compared to the relation between the odometer and the left shift.

 A {\em $\mathbb Z^D$--cover of $M_0$} is a surface $M$ together with a continuous  map $p:M\to M_0$ such that $p$ is a local isometry at every point,  the group of {\em deck transformations}
$$
G:=\{A:M\to M|D\textrm{ an isometry s.t. }p\circ A=p\}
$$
is isomorphic to $\mathbb Z^D$, and for every $x\in M_0$, $p^{-1}(x)$ is a $G$--orbit of some point in $M$.

The flows $g,h:T^1 M_0\to T^1 M_0$ lift to flows $g,h:T^1 M\to T^1 M$ which commute with the elements of $G$, and which satisfy (a),(b). Now (a) and (b) could be compared to the relation between the HIK transformation and a $\mathbb Z^D$--skew-product over the left shift map \cite{Po}.

The locally finite ergodic invariant measures for $h$ are described in \cite{BL} and \cite{S}. There are infinitely many, but only one  up to normalization, is non-squashable \cite{LS}. This measure, which we call $m_0$, is rationally ergodic, and it is  invariant under the action of the geodesic flow and the deck transformations.

We choose a normalization for $m_0$ as follows. Let $\widetilde{M}_0$ be a connected pre-compact subset of $M$ s.t. $p:\widetilde{M_0}\to M_0$ is one-to-one and onto, then we normalize $m_0$ so that $m_0[T^1 \widetilde{M}_0]=1$.

The following can be extracted from \cite{LS}:

\newpage
\proclaim{Theorem 6.1}
\

There exists $a(T)\propto T/(\ln T)^{D/2}$ such that  for every $f\in L^1(m_0)$ with positive integral,
$$
\frac{1}{a(T)}\int_0^T\!\!\! f[h^s(\om)] ds\xrightarrow[T\to\infty]{\frak d}(2^{\frac{D}{2}}e^{-\frac{1}{2}\chi_D^2})m_0(f).
$$
\endproclaim
\noindent

\demo{Proof sketch}
 \
 
 Enumerate $G=\{A_{\underline{\xi}}:\underline{\xi}\in\mathbb Z^D\}$ such that  $A_{\underline{\xi}_1}\circ A_{\underline{\xi}_2}=A_{\underline{\xi}_1+\underline{\xi}_2}$, then $M=\bigcupdot_{\underline{\xi}\in\mathbb Z^D} A_{\underline{\xi}}[\widetilde{M}_0]$.
The {\em $\mathbb Z^D$--coordinate} of $\vec{v}\in T^1 M$ is the unique $\underline{\xi}(\vec{v})\in\mathbb Z^D$ such that  $\vec{v}\in T^1[ A_{\underline{\xi}}(\widetilde{M}_0)]$.

It is known that $\frac{1}{\sqrt{T}}{\underline{\xi}\circ g^T}\xrightarrow[T\to\infty]{\frak d}\mathcal N$, where $\mathcal N$ is a   $D$--dimensional Gaussian random variable  with positive definite covariance matrix $\mathrm{Cov}(\mathcal N)$ (Ratner \cite{R}, Katsuda \& Sunada \cite{KS}).

Let $\|\cdot\|_H$ denote the norm on $\mathbb R^D$ given by $\|\un{v}\|_H:=\sqrt{\underline{v}^t\mathrm{Cov}(\mathcal N)^{-1}\underline{v}}$. The following is proved in \cite{LS} (Theorem 5): Suppose $f\in L^1(m_0)$, then for every $\e>0$, for $m_0$--a.e. $\vec{v}\in T^1 M$, for all $T$ large enough
$$
2^{\frac{D}{2}-\e}e^{-\frac{1}{2}(1+\e)\left\|\frac{\underline{\xi}(g^{\log_\mu T}\vec{v})}{\sqrt{\log_\mu T}}\right\|_H^2}\leq\frac{1}{a(T)}\int_0^T\!\!\!\! f[h^s(\vec{v})]ds\leq
2^{\frac{D}{2}+\e}e^{-\frac{1}{2}(1-\e)\left\|\frac{\underline{\xi}(g^{\log_\mu T}\vec{v})}{\sqrt{\log_\mu T}}\right\|_H^2}
$$
where $a(T)=\const T/(\ln T)^{D/2}$ (the value of the constant is known, see \cite{LS}).

This is the version of   (\Football) needed  to deduce the theorem as above.

\end{document}